\documentclass[a4paper,10pt]{article}
\usepackage{amsmath,amssymb}
\usepackage[utf8]{inputenc}
\usepackage[english,frenchb]{babel}
\usepackage{graphicx}

\newcommand{\ass}{\operatorname{Ass}}
\newcommand{\comm}{\operatorname{Comm}}
\newcommand{\lie}{\operatorname{Lie}}
\newcommand{\mld}{\operatorname{Mould}}
\newcommand{\arbu}{\operatorname{Arb}}
\newcommand{\nap}{\operatorname{NAP}}
\newcommand{\zinb}{\operatorname{Zinb}}

\newcommand{\sym}{\mathfrak{S}}
\newcommand{\NN}{\mathbb{N}}
\newcommand{\ZZ}{\mathbb{Z}}
\newcommand{\ram}{\operatorname{Ram}}
\newcommand{\tfun}{\mathsf{t}}
\newcommand{\ot}{\leftarrow}
\newcommand{\tri}{\triangleleft}

\newcommand{\mex}[1]{\includegraphics[height=5mm]{M#1.eps}}

\newcommand{\ideal}[2]{\mathsf{J}_{#1}(#2)}
\newcommand{\som}[1]{u[#1]}

\newtheorem{theorem}{Théorème}[section] 
\newtheorem{proposition}[theorem]{Proposition}

\newtheorem{lemma}[theorem]{Lemme}

\newenvironment{preuve}{\begin{trivlist}\item{\bf{Preuve.}}}
  {\hfill\rule{2mm}{2mm}\end{trivlist}}

\title{Une opérade anticyclique sur les arbustes}
\author{F. Chapoton}
\date{\today}

\begin{document}
\maketitle

\begin{abstract}
  On définit des objets combinatoires nouveaux, les arbustes, qui
  généralisent les forêts d'arbres enracinés. On introduit ensuite une
  structure d'opérade sur les arbustes. On montre que cette opérade
  $\arbu$ s'injecte dans l'opérade Zinbiel, en utilisant l'inclusion
  de Zinbiel dans l'opérade des moules. On montre aussi que cette
  inclusion est compatible avec la structure plus riche d'opérade
  anticyclique qui existe sur Zinbiel et sur les moules.
\end{abstract}

\selectlanguage{english}

 \begin{abstract}
   We define new combinatorial objects, called shrubs, such that
   forests of rooted trees are shrubs. We then introduce a structure
   of operad on shrubs. We show that this operad is contained in the
   Zinbiel operad, by using the inclusion of Zinbiel in the operad of
   moulds. We also prove that this inclusion is compatible with the
   richer structure of anticyclic operad that exists on Zinbiel and on
   moulds.
\end{abstract}

\selectlanguage{french}

\section*{Introduction}

Parmi les opérades classiques, on trouve l'opérade $\comm$ des
algèbres commutatives et associatives. Cette opérade s'injecte dans
une opérade d'introduction plus récente, l'opérade $\zinb$ des
algèbres de Zinbiel, introduite par Loday \cite{loday_lecture_notes}.
L'opérade $\zinb$ est munie d'une structure plus riche d'opérade
anticyclique \cite{anticyclic}. L'image de $\comm$ dans $\zinb$ n'est
pas une sous-opérade anticyclique et la plus petite sous-opérade
anticyclique de $\zinb$ contenant $\comm$ est en fait $\zinb$ tout
entière.

On peut cependant se poser la question suivante. L'opérade $\comm$ est
en fait une opérade ensembliste, \textit{i.e.} qui existe dans la
catégorie des ensembles plutôt que dans celle des espaces vectoriels.
En oubliant la structure additive de $\zinb$, on obtient une inclusion
d'opérades ensemblistes de $\comm$ dans $\zinb$. Quelle est alors la
plus petite opérade anticyclique ensembliste de $\zinb$ contenant
$\comm$ ? Cet article donne la réponse à cette question, en décrivant
explicitement cette opérade.

On a choisi de présenter les choses dans le sens contraire de la
motivation qui vient d'être donnée. On introduit donc une opérade
\textit{a priori} ; on étudie cette opérade et on montre en particulier qu'elle
s'injecte bien dans $\zinb$ et qu'elle répond à la question posée.

L'opérade en question, qu'on notera $\arbu$, est une opérade
ensembliste, dont la combinatoire fait intervenir des objets nouveaux,
qui généralisent les forêts d'arbres enracinés. On a choisi de nommer
ces objets des arbustes. La première section est consacrée à une étude
combinatoire de ces objets, qui a par ailleurs fait
l'objet d'une note \cite{chapoton08}.

Ensuite, on introduit la structure d'opérade, et on en donne une
présentation par générateurs et relations. L'opérade $\arbu$ est
binaire et quadratique, ce qui signifie qu'elle est engendrée par des
éléments en degré $2$ et que les relations sont en degré $3$.

On décrit ensuite le morphisme d'opérade de $\arbu$ dans $\zinb$ dont
l'existence résulte immédiatement de la présentation par générateurs
et relations.

Pour montrer que ce morphisme est bien injectif, on va utiliser une
injection de $\zinb$ dans une opérade plus grosse, l'opérade des
moules, qu'on notera $\mld$, et dont les éléments sont des fractions
rationnelles en plusieurs variables \cite{moules}.

On donne une description directe du morphisme composé de $\arbu$ dans
$\mld$ en associant explicitement à chaque arbuste une fraction
rationnelle en des variables indexées par ses sommets.

En utilisant alors les descriptions explicites obtenues, on montre que ces
morphismes d'opérades sont injectifs. 

Il reste ensuite à introduire la structure plus riche d'opérade
anticyclique. On introduit pour cela des signes, et on montre que la version
signée de $\arbu$ est une sous opérade anticyclique ensembliste de
$\zinb$. Il est alors évident de voir que c'est bien la plus petite
sous-opérade anticyclique ensembliste de $\zinb$ qui contienne
$\comm$. Ceci répond à la question posée plus haut, qui est la
motivation originelle de cet article.

On termine en donnant la présentation binaire quadratique de l'opérade
duale de $\arbu$, sans l'étudier plus avant.

\section{Combinatoire des arbustes}

Certains aspects purement combinatoires des arbustes ont été présentés
auparavant dans la note \cite{chapoton08}. Par souci de complétude, on
rappelle dans cette section ceux qui seront nécessaires, ainsi que
leurs démonstrations.

\subsection{Définitions}

Un \textbf{arbuste} $P$ sur un ensemble fini $I$ est la donnée d'un
ensemble d'arêtes $i-j$ et d'une fonction hauteur $h_P : I \to \NN$
satisfaisant les conditions suivantes :
\begin{enumerate}
\item Si $i-j$ est une arête, alors $h_P(i)=h_P(j)\pm 1$.
\item Si $h_P(j)>0$, alors il existe une arête $i-j$ avec $h_P(i)=h_P(j)-1$.
\item Les motifs $\mex{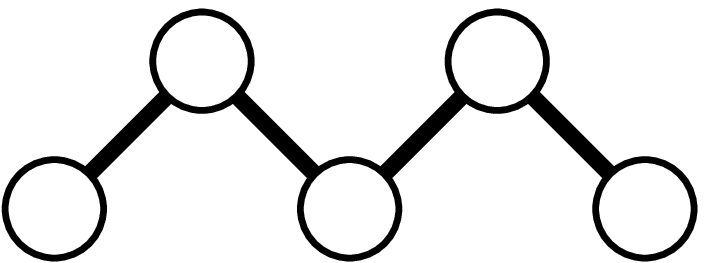}$ et $\mex{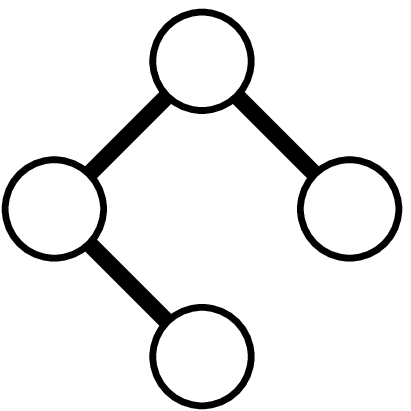}$ sont interdits.
\end{enumerate}

Par motif interdit, on entend que le graphe induit sur un
sous-ensemble de sommets n'est jamais identique à la configuration
exclue, modulo un décalage éventuel de la hauteur.

\textbf{Nota Bene} : on dessine les arbustes avec la hauteur
croissante de bas en haut. En particulier, ceci vaut pour les deux
motifs exclus.

\begin{figure}
  \begin{center}
    \includegraphics[scale=0.35]{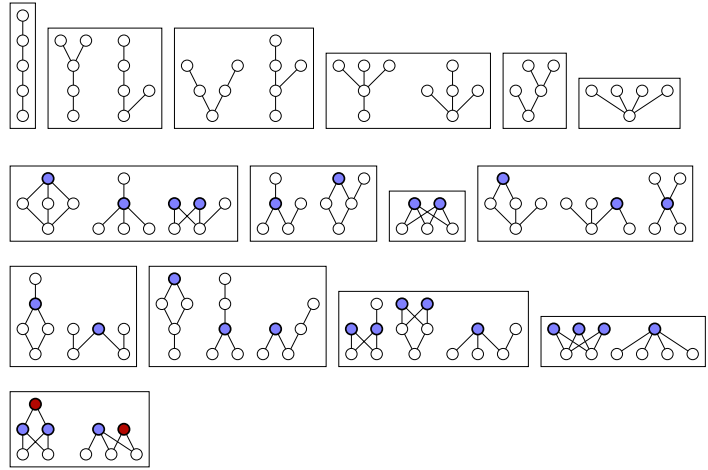} 
    \caption{Les $30$ arbustes connexes sur $5$ sommets.}
    \label{arbuste5}
  \end{center}
\end{figure}

\textbf{Remarque} : la troncation d'un arbuste aux hauteurs plus grandes qu'une
hauteur fixée est encore un arbuste, après décalage de la fonction hauteur.

Les arbres enracinés (\textit{i.e.} les graphes connexes
et simplement connexes munis d'un sommet distingué appelé la racine)
sont des exemples d'arbustes, pour la fonction hauteur donnée par la distance à
la racine. Plus généralement, les forêts d'arbres enracinés sont aussi des
arbustes, pour la même fonction hauteur.

Les graphes bipartites complets, formés d'un ensemble
$I \sqcup J$ et de toutes les arêtes possibles entre $I$ et $J$,
sont aussi des arbustes, en plaçant $I$ à hauteur $0$ et $J$ à hauteur $1$.

\medskip

Par convention, on orientera les arêtes selon la hauteur décroissante.
Si $i-j$ est une arête avec $h_P(i)=h_P(j)-1$, on dira que $j$ couvre
$i$, ce qu'on notera $i \stackrel{P}{\ot} j$ ou simplement $i {\ot} j$ .

Un arbuste est dit \textbf{connexe} si le graphe sous-jacent est
connexe. Tout arbuste s'écrit de manière unique comme union disjointe
d'arbustes connexes.

On appelle \textbf{arbuste trivial} un arbuste à un seul
élément.

La figure \ref{arbuste5} montre les $30$ classes d'isomorphisme
d'arbustes connexes sur $5$ sommets. La signification des couleurs et
des regroupements dans cette figure sera expliquée plus loin.

\medskip

Un sommet $i$ d'un arbuste $P$ est dit \textbf{ramifié} si il est
source d'au moins deux arêtes sortantes, \textit{i.e.} si il couvre au
moins deux sommets distincts. Dans la figure \ref{arbuste5}, les
sommets ramifiés sont coloriés.

Les forêts d'arbres enracinés sont exactement les arbustes sans sommet ramifié.

\subsection{Quelques lemmes}

\begin{lemma}
  \label{chemin_court}
  Soit $P$ un arbuste. Soient $a$ et $a'$ deux éléments de même
  hauteur de $P$. Si il existe un élément $b$ et deux chemins
   \begin{equation}
     a \ot \dots \ot b \to \dots \to a',
   \end{equation}
   alors il existe un élément $b'$ et deux arêtes
   \begin{equation}
     a \ot b'\to a'.
   \end{equation}
\end{lemma}

\begin{preuve}
  Par récurrence sur la hauteur de $b$. Si $b$ est de hauteur
  $h_P(a)+1$, il n'y a rien à démontrer. Sinon, on a des
  élément $c,d,e$ et deux chemins
   \begin{equation}
     a \ot \dots \ot c \ot b \to d \to e \to \dots \to a',
   \end{equation}
   où $e$ peut éventuellement se confondre avec $a'$. On utilise
   alors le motif exclu $\mex{4}$ sur les sommets $\{b,c,d,e\}$ pour
   montrer qu'il existe une arête $c \to e$, donc un chemin de $c$
   à $a'$. Comme la hauteur de $c$ est strictement inférieure à
   celle de $b$, on conclut ensuite par hypothèse de récurrence.
\end{preuve}

\begin{lemma}
  \label{zigzag_court}
  Soit $P$ un arbuste. Soient $a$ et $a'$ deux éléments de
  hauteur nulle de $P$. On suppose qu'il existe des éléments de
  hauteur nulle $a=a_0,\dots,a_\ell=a'$, des éléments $b_i$ de
  hauteur $1$ pour $0 \leq i \leq \ell-1$ et des arêtes
   \begin{equation}
     a_i \ot b_i \to a_{i+1}.
   \end{equation}
   Alors il existe un élément $c$ de hauteur $1$ et des arêtes
   \begin{equation}
     a \ot c \to a'.
   \end{equation}
\end{lemma}
\begin{preuve}
  Par récurrence sur $\ell$. Si $\ell=1$, il n'y a rien à
  démontrer. Sinon on utilise le motif exclu $\mex{5}$ sur les
  éléments $a_0,a_1,a_2,b_0,b_1$. On en déduit soit qu'il
  existe une arête de $b_0$ vers $a_2$, soit qu'il existe une
  arête de $b_1$ vers $a_0$. 

  Dans tous les cas, on obtient, en supprimant $a_1$ et $b_0$ ou
  $b_1$, une suite plus courte satisfaisant les mêmes
  hypothèses. On conclut par récurrence.
\end{preuve}

\begin{lemma}
  \label{paire_zero}
  Soit $P$ un arbuste connexe et non trivial sur l'ensemble $I$.
  Soient $a$ et $a'$ distincts de hauteur nulle dans $P$. Alors il
  existe $b$ couvrant $a$ et $a'$.
 \end{lemma}

 \begin{preuve}
   Soient $a$ et $a'$ de hauteur nulle. Par connexité de $P$, on peut
   trouver une suite d'arêtes reliant $a$ et $a'$. Comme tout élément
   de hauteur non nulle est source d'au moins une arête sortante, on
   peut, quitte à prolonger cette suite d'arêtes, trouver des éléments
   de hauteur nulle $a=a_0,\dots,a_\ell=a'$, des éléments $b_i$ pour
   $0 \leq i \leq \ell-1$ et des chemins
   \begin{equation}
     a_i \ot \dots \ot b_i \to \dots \to a_{i+1}.
   \end{equation}
   
   Par le lemme \ref{chemin_court}, on peut trouver des éléments
   $b'_i$ de hauteur $1$ tels que
   \begin{equation}
     a_i  \ot b'_i \to a_{i+1}.
   \end{equation}

   Par le lemme \ref{zigzag_court}, il existe $b$ couvrant $a$ et
   $a'$.
 \end{preuve}

\begin{proposition}
  \label{existe_couvrant}
  Soit $P$ un arbuste connexe et non trivial sur l'ensemble $I$. Il
  existe dans $P$ au moins un élément de hauteur $1$ lié à
  tous les éléments de hauteur nulle.
\end{proposition}

\begin{preuve}
  Prenons $e$ de hauteur $1$ couvrant le plus possible
  d'éléments de hauteur nulle. Supposons qu'il existe $a$ de
  hauteur nulle qui ne soit pas couvert par $e$.

  Soit $b$ un élément couvert par $e$. Soit $d$ une couverture commune
  à $a$ et $b$, qui existe par le lemme \ref{paire_zero}. Alors
  $d$ diffère de $e$ par hypothèse.

  Soit $c$ un élément couvert par $e$ et différent de $b$, si
  il en existe. Considérons le graphe induit sur $\{a,b,c,d,e\}$.
  Par le motif exclu $\mex{5}$, ou bien $e$ couvre $a$, ou bien $d$
  couvre $c$. Comme $e$ ne couvre pas $a$ par hypothèse, on a donc
  que $d$ couvre $c$.

  Qu'un tel $c$ existe ou non, on a donc montré que $d$ couvre tous
  les éléments couverts par $e$ et aussi $a$, ce qui contredit
  la maximalité de $e$. Donc l'existence de $a$ est absurde et $e$
  couvre tous les éléments de hauteur nulle.
\end{preuve}

\begin{proposition}
  \label{deconnexion}
  Soit $P$ un arbuste connexe et non trivial sur l'ensemble $I$. Soit
  $S$ l'ensemble des éléments de hauteur $1$ couvrant tous les
  éléments de hauteur nulle. Dans le graphe obtenu en ôtant dans $P$
  toutes les arêtes entre $S$ et les éléments de hauteur nulle, aucun
  sommet de $S$ n'est dans la même composante connexe qu'un sommet de
  hauteur nulle.
\end{proposition}

\begin{preuve}
  Supposons le contraire. Soit $P_{\geq 1}$ l'arbuste obtenu par
  troncation de $P$ aux hauteurs non nulles. Il existe donc un élément
  $a$ de $S$ dont la composante connexe dans $P_{\geq 1}$ contient un
  élément $a'$ de $h_P^{-1}(1) \setminus S$. Par le lemme
  \ref{paire_zero} appliqué à cette composante connexe, il existe $b$
  et des arêtes $a \ot b \to a'$.

  Par le motif exclu $\mex{4}$ appliqué dans $P$ à $a,b,a'$ et un
  élément $c$ de hauteur nulle, on montre alors que $a'$ couvre tous
  les éléments de hauteur nulle de $P$, donc appartient à $S$, ce qui
  est absurde.
\end{preuve}

On appelle \textbf{feuille} un sommet ayant une seule arête sortante
et aucune arête entrante. On dit que deux sommets sont
\textbf{corrélés} si ils ont le même ensemble de sources des arêtes
entrantes et le même ensemble de buts des arêtes sortantes.

\begin{lemma}
  \label{feuille_ou_correles}
  Dans un arbuste non-trivial, il existe au moins une feuille ou deux
  sommets corrélés.
\end{lemma}

\begin{preuve}
  Supposons qu'il n'existe pas de feuille et montrons qu'il existe
  deux sommets corrélés. En effet, soit $c$ un sommet de ramification
  minimale parmi les sommets de hauteur maximale. Par hypothèse, $c$
  n'est pas une feuille, donc est ramifié. Soient $a$ et $b$ deux
  sommets distincts couverts par $c$. Montrons que $a$ et $b$ sont
  corrélés.

  Soit $d$ un sommet de même hauteur que $c$ qui couvre soit $a$ soit
  $b$. Montrons que $d$ couvre $a$ et $b$.  Supposons par exemple que
  $d$ ne couvre pas $a$ mais seulement $b$. Par minimalité de la
  ramification de $c$, il existe un sommet $e$ qui est couvert par $d$
  mais pas couvert par $c$. On en déduit le motif $\mex{5}$ sur
  $\{c,d,a,b,e\}$, ce qui est absurde.

  Par ailleurs, tout sommet $e$ couvert par $a$ ou $b$ est en fait
  couvert par $a$ et $b$, en raison du motif exclu $\mex{4}$ appliqué
  aux sommets $\{c,a,b,e\}$.
\end{preuve}

\begin{lemma}
  Soit $P$ un arbuste et $i$ une feuille de $P$. Le graphe obtenu en
  supprimant le sommet $i$ est un arbuste.
\end{lemma}
\begin{preuve}
  C'est clair. La fonction hauteur est inchangée sur les autres
  sommets. Les motifs exclus le sont toujours.
\end{preuve}

\begin{lemma}
  Soit $P$ un arbuste et $i,j$ deux sommets corrélés de $P$. Le
  graphe obtenu en identifiant $i$ et $j$ est un arbuste.
\end{lemma}
\begin{preuve}
  La fonction hauteur passe au quotient. Si un des motifs exclus
  apparaissait dans le graphe quotient, on pourrait le relever dans $P$ en
  choisissant de relever en $i$ le sommet obtenu par identification.
\end{preuve}

\subsection{Bijection avec les posets série-parallèle}

On peut construire par récurrence une bijection entre les arbustes et
les posets Série-Parallèle, voir \cite{stanley1974,stanleyEC2} pour
cette notion. Plus précisément, il existe un isomorphisme entre
l'espèce des arbustes et celle des posets Série-Parallèle, au sens de
la théorie des espèces de structure de Joyal \cite{Bible}. On renvoie
le lecteur à la note \cite{chapoton08} pour les détails de cette
construction.

On peut noter qu'il existe aussi une structure d'opérade sur les
posets Série-Parallèle, voir l'article \cite{loday-2008} de Loday et
Ronco.





\section{L'opérade des arbustes}

\subsection{Définition}

Dans cette section, on définit une structure d'opérade en se
donnant des compositions partielles $\circ_i$ sur les arbustes.

Soit $P$ un arbuste sur $I$, $P'$ un arbuste sur $I'$ et $i \in I$.

Définissons d'abord une fonction hauteur sur $I \setminus \{i\}
\sqcup I'$ :
\begin{equation}
  h_{P \circ_i P'} (j)=
  \begin{cases}
    h_P(j) &\text{ si } j\in I \setminus \{i\},\\
    h_{P'}(j)+h_P(i) &\text{ si } j\in I'.
  \end{cases}
\end{equation}

\begin{figure}
  \begin{center}
    \includegraphics[scale=0.35]{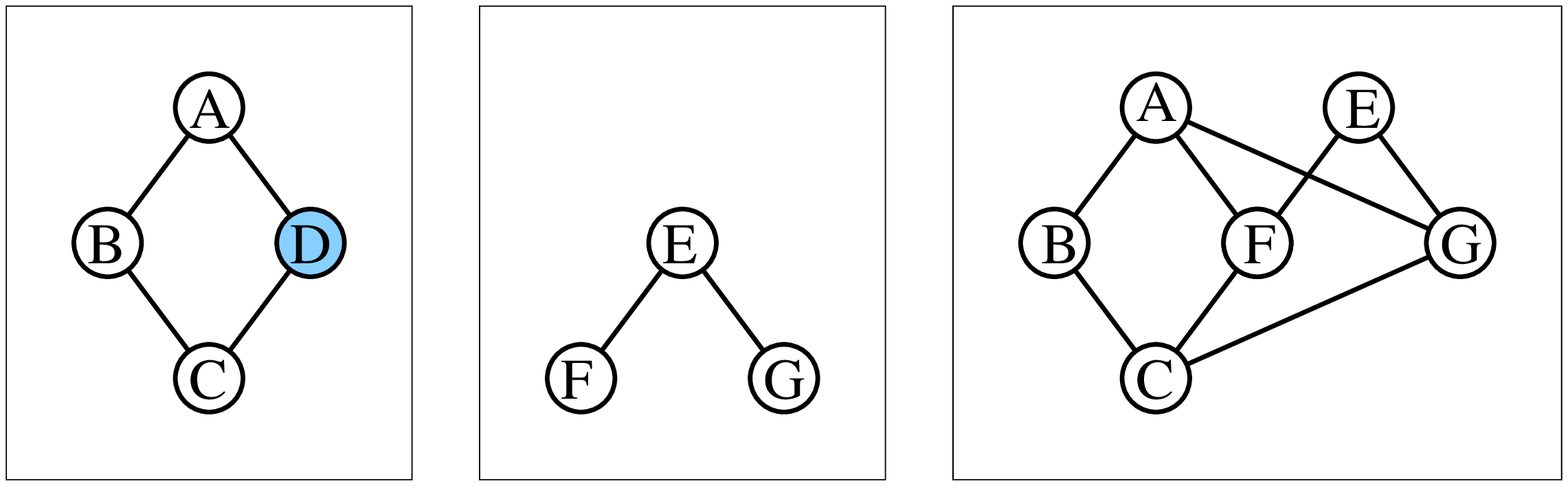} 
    \caption{Composition de deux arbustes $P$ et $Q$ : $P \circ_{D} Q$.}
    \label{composi}
  \end{center}
\end{figure}

Définissons aussi un graphe $P \circ_i P'$ sur $I \setminus \{i\}
\sqcup I'$ comme suit : dans l'union disjointe de $P'$ et de $P$ privé
du sommet $i$, on relie par une arête tout sommet de $P$ qui était
relié à $i$ à tout sommet de hauteur nulle de $P'$. La figure
\ref{composi} montre un exemple de cette construction.

\begin{proposition}
  Le graphe $P \circ_i P'$ (muni de la hauteur $h_{P \circ_i P'}$) est un arbuste.
\end{proposition}

\begin{preuve}
  Les conditions sur la hauteur sont clairement vérifiées.  Il reste
  donc à vérifier que les motifs exclus $\mex{4}$ et $\mex{5}$
  n'apparaissent pas.

  \textbf{Remarque préliminaire} : supposons qu'un des deux motifs
  exclus soit présent dans $P \circ_i P'$. Nécessairement, comme le
  motif est exclu de $P$ et de $P'$, une des arêtes du motif relie un
  élément de $I\setminus \{i\}$ à un élément de $I'$. Si un seul
  sommet du motif est dans $I'$, alors on peut relever le motif dans
  $P$ en remplaçant ce sommet par $i$, ce qui est absurde. Donc, au
  moins deux sommets du motif sont dans $I'$.

  Supposons d'abord qu'on ait le motif $\mex{5}$ avec $a,b$ en haut et
  $c,d,e$ en bas. Par la remarque ci-dessus, la hauteur de $a$ et $b$
  est donc soit $h_P(i)$ soit $h_P(i)+1$.

  Supposons d'abord que la hauteur de $a$ et $b$ est $h_P(i)$. Alors
  $c,d,e$ sont dans $I\setminus\{i\}$. Par la remarque ci-dessus, $a$
  et $b$ sont dans $I'$. Mais dans ce cas, les cinq sommets
  $a,b,c,d,e$ ne peuvent pas former le motif $\mex{5}$, par
  construction de la composition.
 
  Supposons ensuite que la hauteur de $a$ et $b$ est $h_P(i)+1$. Si
  $a$ est dans $I'$, alors $c,d$ sont dans $I'$ par construction. Mais
  dans ce cas, $b$ est aussi dans $I'$, sinon, il y aurait une arête
  $b-c$ par construction. Par conséquent $a,b,c,d,e$ sont dans $I'$ et
  $P'$ contient le motif $\mex{5}$, ce qui est absurde. On montre de
  même que $b\not\in I'$. On peut donc supposer que $a$ et $b$ sont
  dans $I\setminus\{i\}$. Nécessairement, on a alors que $c$ et $e$
  sont dans $I\setminus\{i\}$. On aurait alors au plus un sommet dans
  $I'$, ce qui est impossible par la remarque.

  \medskip

  Supposons maintenant qu'on ait le motif $\mex{4}$ avec $a$ en haut,
  $b,c$ au milieu et $d$ en bas. Par la remarque, la hauteur de $a$
  est $h_P(i),h_P(i)+1$ ou $h_P(i)+2$.

  Si $a \in I'$ est de hauteur $h_P(i)$, alors $b,c,d$ sont dans
  $I\setminus\{i\}$, ce qui contredit la remarque.

  Si $a \in I'$ est de hauteur $h_P(i)+1$, alors $b,c$ sont dans $I'$
  et $d$ est dans $I\setminus\{i\}$. Le motif est impossible par
  construction de la composition.

  Si $a \in I'$ est de hauteur $h_P(i)+2$, alors $b,c,d$ sont dans
  $I'$, ce qui est absurde car le motif $\mex{4}$ est exclu dans $P'$.

  On a donc nécessairement $a\in I\setminus\{i\}$. Alors, par la
  remarque, la hauteur de $a$ est $h_P(i)+1$, $b,c$ sont dans $I'$ et
  $d$ dans $I\setminus\{i\}$. Mais alors le motif $\mex{4}$ est
  impossible par construction de la composition.

\end{preuve}

\begin{proposition}
  \label{operade_narb}
  Ces compositions définissent une opérade $\arbu$.
\end{proposition}

\begin{preuve}
  L'arbuste trivial est clairement une unité pour la composition.

  La composition est aussi clairement un application fonctorielle, qui
  ne dépend en aucune façon des étiquettes des sommets.

  Il reste à vérifier les deux axiomes d'associativité des
  opérades.

  Soient $P,P'$ et $P''$ des arbustes sur $I,I'$ et $I''$. Soient $i$ et
  $j$ dans $I$. Si $i$ et $j$ ne sont pas reliés dans $P$, alors on
  a clairement
  \begin{equation}
    (P \circ_i P') \circ_j P'' = (P \circ_j P'') \circ_i P',
  \end{equation}
  dans la mesure où la composition est locale dans $P$.

  Si $i$ et $j$ sont reliés par une arête, cette relation est
  encore vraie, le résultat contenant un graphe bipartite complet
  entre les sommets de hauteur nulle de $P'$ et ceux de $P''$.

  Soient $P,P'$ et $P''$ des arbustes sur $I,I'$ et $I''$. Soient $i$
  dans $I$ et $i'$ dans $I'$. Alors on a bien
  \begin{equation}
    (P \circ_i P') \circ_{i'} P'' = P \circ_i (P' \circ_{i'} P''). 
  \end{equation}
  En effet, la description du résultat est la même pour les deux
  expressions, en distinguant deux cas selon que $i'$ est de hauteur
  nulle ou non dans $P'$. 

  Si la hauteur de $i'$ dans $P'$ est nulle, on relie tous les sommets
  de hauteur nulle de $P''$ à tous les sommets reliés à $i$ dans $P$
  et à tous les sommets reliés à $i'$ dans $P'$ et on relie tous les
  sommets de hauteur nulle de $P'$ à tous les sommets liés à $i$ dans
  $P$.

  Si la hauteur de $i'$ dans $P'$ est non nulle, on relie tous les
  sommets de hauteur nulle de $P''$ à tous les sommets liés à $i'$
  dans $P'$ et on relie tous les sommets de hauteur nulle de $P'$ à
  tous les sommets liés à $i$ dans $P$.
\end{preuve}

On note $[1][2]$ l'union disjointe de deux arbustes triviaux.
L'opération correspondante $(P,Q)\mapsto [1][2] \circ_1 P \circ_2 Q$
est juste l'union disjointe des arbustes. C'est un produit associatif
et commutatif. On notera simplement $P \, Q$ l'union disjointe de $P$
et $Q$.

On note $[1\tri 2]$ l'arbre enraciné sur l'ensemble $\{1,2\}$ ayant
pour racine $1$. L'opération correspondante $(P,Q)\mapsto [1 \tri 2]
\circ_1 P \circ_2 Q$ est la suivante : on joint les éléments de
hauteur $0$ de $P$ à ceux de hauteur $0$ de $Q$ par un graphe
bipartite complet, et on ajoute $1$ à la fonction hauteur de $Q$.
Cette opération n'est pas associative. On la notera $P \tri Q$.

\textbf{Remarque} : on a deux sous-opérades simples contenues dans $\arbu$.
Les arbres enracinés forment une copie de l'opérade $\nap$
\cite{NAP_muriel}. Les forêts d'arbustes triviaux forment une copie
de l'opérade $\comm$.

\subsection{Présentation par générateurs et relations}

\begin{theorem}
  \label{presentation}
  L'opérade $\arbu$ est engendrée par ses éléments $[1][2]$ et $[1 \tri 2]$ de
  degré $2$, soumis uniquement aux relations suivantes :
  \begin{equation}
    \label{relaNAP}
    [* \tri 1] \circ_* [3 \tri 2]=[* \tri 2] \circ_* [3 \tri 1]=[3 \tri *] \circ_* [1][2]
  \end{equation}
  et
  \begin{equation}
    \label{relaComm}
    [*][1] \circ_* [2][3]= [*][2] \circ_* [3][1].
  \end{equation}
\end{theorem}

\begin{preuve}
  Soit $\mathtt{A}$ l'opérade quotient de l'opérade libre sur les
  générateurs $C_{1,2}=C_{2,1}$ et $D_{1,2}$ en degré $2$ par l'idéal
  engendré par les relations obtenues en remplaçant $[1][2]$ par
  $C_{1,2}$ et $[1 \tri 2]$ par $D_{1,2}$ dans les relations de
  l'énoncé.

  Tout d'abord, les relations de l'énoncé sont faciles à vérifier dans
  $\arbu$. On a donc un unique morphisme $\phi$ de l'opérade
  $\mathtt{A}$ dans $\arbu$ tel que $\phi(C_{1,2})=[1][2]$ et
  $\phi(D_{1,2})=[1 \tri 2]$. Il faut montrer que ce morphisme est un
  isomorphisme en construisant un inverse $\psi$.

  On construit $\psi$ par récurrence sur le degré $n$. On suppose par
  récurrence que $\psi$ est l'inverse de $\phi$ jusqu'en degré $n-1$.
  Ceci entraîne en particulier que $\psi$ est un morphisme d'opérade
  jusqu'en degré $n-1$.

  Soit $P$ un arbuste de degré $n\geq 2$. Par le Lemme
  \ref{feuille_ou_correles}, il existe dans $P$ une feuille ou une
  paire de sommets corrélés.

  Si il existe une feuille $i$, on définit $\psi$ de la façon
  suivante. Soit $j$ l'unique sommet couvert par $i$. Soit $P'$
  l'arbuste obtenu à partir de $P$ en supprimant le sommet $i$ et en
  ré-étiquetant par le symbole $\square$ le sommet $j$. Alors $P$
  peut s'écrire $P' \circ_{\square} [j \tri i]$. On pose
  $\psi(P)=\psi(P')\circ_{\square} D_{j,i}$.
  
  Si il existe deux sommets corrélés $i$ et $j$, on définit $\psi$ de
  la façon suivante. Soit $P'$ l'arbuste obtenu à partir de $P$ en
  identifiant les sommets $i$ et $j$ de $P$ et en étiquetant le
  sommet ainsi défini par $\square$. L'arbuste $P$ peut alors s'écrire
  $P' \circ_{\square} [i][j]$. On pose
  $\psi(P)=\psi(P')\circ_{\square} C_{i,j}$.

  Il se peut bien entendu qu'il existe plusieurs feuilles, plusieurs
  paires de sommets corrélés ou qu'il existe à la fois une feuille et
  une paire de sommets corrélés dans $P$. Il faut montrer que la
  définition de $\psi$ est bien cohérente, \textit{i.e.} que nous
  avons bien ainsi défini une application $\psi$, indépendamment de
  tout choix. On va considérer successivement le cas de deux feuilles,
  le cas de deux paires et le cas d'une feuille et d'une paire, et
  deux sous-cas pour chacun de ces trois cas.

  \medskip

  1) Soient $i$ et $i'$ deux feuilles. Si elles couvrent deux sommets
  distincts $j$ et $j'$, alors il existe un arbuste $P''$ tel que
  \begin{equation}
    P = (P'' \circ_{*} [j' \tri i']) \circ_{\square} [j \tri i] = (P'' \circ_{\square} [j \tri i]) \circ_{*} [j' \tri i'].
  \end{equation}
  On obtient $P''$ en supprimant $i$ et $i'$ dans $P$ et en
  ré-étiquetant $j$ par $\square$ et $j'$ par $*$.

  Dans ce cas, les deux définitions possibles de $\psi$ coïncident :
  \begin{multline}
    \psi(P'' \circ_{*} [j' \tri i'])\circ_{\square} D_{j,i} \\
    = (\psi(P'') \circ_{*} D_{j',i'}) \circ_{\square} D_{j,i} \\ 
    = (\psi(P'') \circ_{\square} D_{j,i}) \circ_{*} D_{j',i'} \\
    =\psi(P''\circ_{\square} [j \tri i])\circ_{*} D_{j',i'},
  \end{multline}
  par hypothèse de récurrence et en utilisant l'axiome de
  commutativité des opérades.

  Si les deux feuilles $i$ et $i'$ couvrent le même sommet $j$, alors
  il existe un arbuste $P''$ tel que
  \begin{equation}
    P = P'' \circ_{\square} ([* \tri i]) \circ_{*} [j \tri i'])
    = P'' \circ_{\square} ([* \tri i'] \circ_{*} [j \tri i]).
  \end{equation}
  On obtient $P''$ en supprimant $i$ et $i'$ dans $P$ et en
  ré-étiquetant $j$ par $\square$.

  Dans ce cas, les deux définitions possibles de $\psi$ coïncident :  
  \begin{multline}
    \psi(P'' \circ_{\square} [* \tri i]) \circ_{*} D_{j , i'}\\
    =(\psi(P'') \circ_{\square} D_{*,i} ) \circ_{*} D_{j,i'}    
    =\psi(P'') \circ_{\square} (D_{*,i}  \circ_{*} D_{j,i'}) \\
    = \psi(P'') \circ_{\square} (D_{*,i'}  \circ_{*} D_{j,i})
    =( \psi(P'') \circ_{\square} D_{*,i'} ) \circ_{*} D_{j,i} \\
    =\psi(P'' \circ_{\square} [* \tri i'] ) \circ_{*} D_{j,i},
  \end{multline}
  par la partie gauche de la relation (\ref{relaNAP}).

  \medskip

  2) Soient $i,j$ et $i',j'$ deux paires de sommets corrélés, tous
  distincts. Alors il existe un arbuste $P''$ tel que
  \begin{equation}
    P = (P'' \circ_{*} [i'][j'] ) \circ_{\square} [i][j] = (P'' \circ_{\square} [i][j]) \circ_{*} [i'][j'].
  \end{equation}
  On obtient $P''$ en identifiant $i$ et $j$ d'une part et $i'$ et
  $j'$ d'autre part dans $P$ et en étiquetant ces sommets par
  $\square$ et $*$.

  Dans ce cas, les deux définitions possibles de $\psi$ coïncident :
  \begin{multline}
    \psi(P'' \circ_{*} [i'][j'])\circ_{\square} C_{i,j}\\
    = (\psi(P'') \circ_{*} C_{i',j'}) \circ_{\square} C_{i,j} \\
    = (\psi(P'') \circ_{\square} C_{i,j}) \circ_{*} C_{i',j'} \\
    =\psi(P''\circ_{\square} [i][j])\circ_{*} C_{i',j'},
  \end{multline}
  par hypothèse de récurrence et en utilisant l'axiome de
  commutativité des opérades.

  Si $i,j$ et $j,k$ sont deux paires de sommets corrélés, alors il
  existe un arbuste $P''$ tel que
  \begin{equation}
    P = P'' \circ_{\square} ([*][i] \circ_{*} [j][k])= P'' \circ_{\square} ([*][k] \circ_{*} [j][i]).
  \end{equation}
  On obtient $P''$ en identifiant les trois sommets $i,j,k$ dans $P$
  et en étiquetant ce sommet par $\square$.

  Dans ce cas, les deux définitions possibles de $\psi$ coïncident :  
  \begin{multline}
    \psi(P'' \circ_{\square} [*] [i]) \circ_{*} C_{j,k}\\
    =(\psi(P'') \circ_{\square} C_{*,i}) \circ_{*} C_{j,k} 
    =\psi(P'') \circ_{\square} (C_{*,i} \circ_{*} C_{j,k}) \\
    =\psi(P'') \circ_{\square} (C_{*,k} \circ_{*} C_{j,i})
   =(\psi(P'') \circ_{\square} C_{*,k}) \circ_{*} C_{j,i} \\
    =\psi(P'' \circ_{\square} [*] [k]) \circ_{*} C_{j,i},    
  \end{multline}
  par la relation (\ref{relaComm}).
  
  \medskip

  3) Soit enfin $i$ une feuille et $j,k$ des sommets corrélés. Si $i$ est
  distincte de $j$ et de $k$, alors $i$ couvre un sommet $\ell$ distinct de
  $j$ et de $k$ et il existe $P''$ tel que
  \begin{equation}
    P = (P'' \circ_{*} [j][k]) \circ_{\square} [\ell \tri i] = (P'' \circ_{\square} [\ell \tri i]) \circ_{*} [j][k].
  \end{equation}
  On obtient $P''$ en supprimant $i$ et en identifiant $j$ et $k$ dans
  $P$, puis en ré-étiquetant $\ell$ par $\square$ et en étiquetant le
  sommet provenant de $j,k$ par $*$.

  Dans ce cas, les deux définitions possibles de $\psi$ coïncident :
  \begin{multline}
    \psi(P'' \circ_{*} [j][k]) \circ_{\square} D_{\ell,i}\\
    = (\psi(P'')\circ_{*} C_{j,k} ) \circ_{\square} D_{\ell,i} \\
    = (\psi(P'') \circ_{\square} D_{\ell,i}) \circ_{*} C_{j,k} \\
    =\psi(P''\circ_{\square} [\ell \tri i] ) \circ_{*} C_{j,k},
  \end{multline}
  par l'axiome de commutativité des opérades. 

  Si $i=j$, alors $k$ est aussi une feuille, les sommets $i$ et $k$
  couvrent un sommet $\ell$ et il existe un arbuste $P''$ tel que
  \begin{equation}
    P = (P'' \circ_{\square} [\ell \tri *]) \circ_{*} [i][k]= (P'' \circ_{\square} [* \tri i] )\circ_{*} [\ell \tri k].
  \end{equation}
  On obtient $P''$ en supprimant $j$ et $k$ dans $P$ et en
  ré-étiquetant le sommet $\ell$ par $\square$.

  Dans ce cas, les deux définitions possibles de $\psi$ coïncident :  
  \begin{multline}
    \psi(P'' \circ_{\square} [\ell \tri *]) \circ_{*} C_{i,k}\\
    =(\psi(P'') \circ_{\square} D_{\ell,*}) \circ_{*} C_{i,k} 
    =\psi(P'') \circ_{\square} (D_{\ell,*} \circ_{*} C_{i,k}) \\
    = \psi(P'') \circ_{\square} (D_{*,i} \circ_{*} D_{\ell,k})
    =( \psi(P'') \circ_{\square} D_{*,i}) \circ_{*} D_{\ell,k}\\
    =\psi(P'' \circ_{\square} [* \tri i]) \circ_{*} [\ell \tri k],    
  \end{multline}
  par la partie droite de la relation (\ref{relaNAP}).

  \medskip

  Ceci termine de montrer que l'on a bien défini une application $\psi$.

  Montrons maintenant que $\psi$ et $\phi$ sont bien inverses l'un de l'autre.

  Comme $\phi$ est un morphisme d'opérade, on a clairement $\phi
  \circ \psi=Id$ par hypothèse de récurrence.

  Considérons une expression monomiale $R$ de degré $n$ dans l'opérade
  $\mathtt{A}$.

  Supposons d'abord qu'elle s'écrive $R=R' \circ_i D_{1,2}$. Alors
  $\phi(R)=\phi(R') \circ_i [1 \tri 2]$. Donc $2$ est une feuille de
  $\phi(R)$, qui couvre $1$. En choisissant cette feuille pour
  calculer $\psi$, on trouve
  \begin{equation}
    \psi(\phi(R))=\psi(\phi(R') \circ_i [1 \tri 2])=\psi(\phi(R')) \circ_i D_{1,2}=R'\circ_i D_{1,2}=R.
  \end{equation}
  
  Sinon, on peut écrire $R=R' \circ_i C_{1,2}$. Alors
  $\phi(R)=\phi(R') \circ_i [1][2]$. Donc $\phi(R)$ possède deux
  sommets corrélés : $1,2$. En choisissant ces sommets corrélés pour
  calculer $\psi$, on trouve
  \begin{equation}
    \psi(\phi(R))=\psi(\phi(R') \circ_i [1][2])=\psi(\phi(R')) \circ_i C_{1,2}=R'\circ_i C_{1,2}=R.
  \end{equation}

  Donc les applications $\phi$ et $\psi$ sont inverses l'une de
  l'autre en degré $n$. Ceci conclut la récurrence.
\end{preuve}

\section{Morphismes injectifs}

\subsection{Morphisme vers l'opérade $\zinb$}

Soit $P$ un arbuste sur un ensemble $I$. On dit qu'un ordre total $<$
sur $I$ est \textbf{compatible} avec $P$ si il satisfait la condition
suivante :
\begin{equation}
  \label{propo}
  \forall i \in I \quad  \exists j \in I  \quad j \stackrel{P}{\ot} i \text{ and } j < i,
\end{equation}
où la notation $j \stackrel{P}{\ot} i$ signifie que $i$ couvre $j$
dans $P$.

Autrement dit, chaque sommet $i$ doit être supérieur pour l'ordre $<$
à au moins un des sommets $j$ qu'il couvre dans $P$.

Si $P$ est non-ramifié, la condition (\ref{propo}) décrit simplement
l'ordre partiel associé à une forêt d'arbres enracinés et les ordres totaux compatibles avec $P$ sont les extensions linéaires de cet ordre
partiel.

\smallskip

On rappelle brièvement la définition de l'opérade $\zinb$, voir
\cite{loday_lecture_notes}. Si $I$ est un ensemble fini, $\zinb(I)$
est l'espace vectoriel de base les ordres totaux sur $I$. On va
considérer un ordre total comme une liste, avec le minimum à gauche.

Soit $\pi\in\zinb(I)$ et $\sigma\in\zinb(I')$ deux listes. On note
$\leq_\pi$ et $\leq_\sigma$ les ordres totaux correspondants sur $I$
et sur $I'$. La composition de $\sigma$ dans $\pi$ en position $i\in
I$ est une somme sur les battages de $\sigma$ avec $\pi$ tels que
$\sigma$ reste à droite de $i$. Plus précisément, soit $\prec$ l'ordre
partiel sur $I\setminus\{i\} \sqcup I'$ défini par
\begin{equation}
  a\preceq b \text{ si }\begin{cases}
    \text{    }&a,b\in I \text{ et }a\leq_\pi b\\
    \text{ ou }&a,b\in I' \text{ et }a\leq_\sigma b\\
    \text{ ou }&a\in I,b\in I'\text{ et }a\leq_\pi i.
  \end{cases}
\end{equation}
Alors la composition $\pi \circ_i \sigma$ dans $\zinb(I\setminus\{i\}
\sqcup I')$ est la somme des ordres totaux qui étendent l'ordre
partiel $\preceq$.

\begin{proposition}
  \label{morphisme_zinbiel}
  Il existe un unique morphisme d'opérades $\gamma$ de $\arbu$ dans
  $\zinb$ tel que 
  \begin{align}
    \gamma( [2 \tri 1] ) & =  [21],\\
    \gamma( [1][2] ) & = [12]+[21].
  \end{align}
\end{proposition}

\begin{preuve}
  Il suffit de vérifier que les relations (\ref{relaNAP}) et
  (\ref{relaComm}) sont satisfaites par les images dans $\zinb$, ce
  qui est immédiat.
\end{preuve}

Le morphisme $\gamma$ admet la description suivante.
\begin{proposition}
  Pour tout arbuste $P$, $\gamma(P)$ est la somme des ordres totaux
  sur $I$ qui sont compatibles avec $P$.
\end{proposition}

\begin{preuve}
  Comme $\arbu$ est engendrée par ses éléments de degré $2$, il suffit
  de montrer par récurrence que c'est vrai pour $[1][2] \circ_1 P
  \circ_2 Q$ et pour $[1 \tri 2] \circ_1 P \circ_2 Q$ si c'est vrai
  pour $P$ et pour $Q$.
  
  L'arbuste $[1][2] \circ_1 P \circ_2 Q$ est simplement l'union
  disjointe de $P$ et $Q$. Les ordres compatibles avec l'union
  disjointe sont exactement les battages d'un ordre $\pi$ compatible
  avec $P$ et d'un ordre $\sigma$ compatible avec $Q$. Or $([12]+[21])
  \circ_1 \pi \circ_2 \sigma$ dans $\zinb$ est exactement la somme des
  battages des deux listes $\pi$ et $\sigma$. On en déduit que la
  proposition est vérifiée pour l'union disjointe de $P$ et $Q$.

  Il reste à considérer le cas de $[1 \tri 2] \circ_1 P \circ_2 Q$.
  C'est l'arbuste obtenu en reliant tous les sommets de hauteur nulle
  de $P$ à tous les sommets de hauteur nulle dans $Q$ et en ajoutant
  $1$ à la fonction hauteur de $Q$. En termes d'ordres compatibles,
  ceci signifie que l'on a les mêmes relations entre les sommets de
  $P$ que pour un ordre compatible avec $P$, les mêmes relations entre
  les sommets de $Q$ que pour un ordre compatible avec $Q$ et la
  condition supplémentaire suivante : pour tout sommet de hauteur
  nulle $i$ de $Q$, il existe un sommet $j$ de hauteur nulle de $P$
  tel que $j < i$. Les ordres compatibles sont donc les battages d'un
  ordre $\pi$ compatible pour $P$ avec un ordre $\sigma$ compatible
  pour $Q$ tels que le minimum est un élément de $P$. Or l'opération
  $[21] \circ_1 \pi \circ_2 \sigma$ dans $\zinb$ est exactement la
  somme des battages des deux listes $\pi$ et $\sigma$ tels que le
  minimum provient de $\pi$.
\end{preuve}

\subsection{Morphisme vers l'opérade des moules}

Commençons par rappeler la définition de l'opérade $\mld$ des moules,
voir \cite{moules,moules2} pour plus de détails.

Si $I$ est un ensemble fini, $\mld(I)$ est le corps des fractions
rationnelles en des indéterminées $(u_i)_{i\in I}$. 

La composition d'une fraction $g$ dans $\mld(J)$ dans une fraction
$f\in\mld(I)$ en position $i\in I$ est donnée par la formule suivante :
\begin{equation}
  f \circ_i g = \big{(}\sum_{j \in J} u_j\big{)}\times g\times f\big{\vert}_{u_i=\sum_{j \in J} u_j} .
\end{equation}
Dans cette formule, $g$ est appliquée aux variables $u_j$ pour $j\in
J$, et on remplace dans les arguments de $f$ la variable $u_i$ par la
somme indiquée.

On a une inclusion de $\zinb$ dans $\mld$ (voir \cite{moules2}) qui
est donnée par l'application suivante : un ordre total $\pi$ sur $I$
est envoyé sur la fraction
\begin{equation}
  \frac{1}{\prod_{i \in I} \sum_{j \geq_\pi i} u_j}.
\end{equation}

\textbf{Remarque} : on voit assez facilement que l'image dans $\mld$ d'un
élément de $\zinb$ permet de calculer cet élément en utilisant des
résidus itérés. Ceci entraîne l'injectivité.

En composant le morphisme $\gamma$ de $\arbu$ dans $\zinb$ avec cette
inclusion de $\zinb$ dans $\mld$, on obtient un morphisme $\kappa$ de $\arbu$
dans $\mld$, qui est caractérisé comme suit.

\begin{proposition}
  \label{morphisme_moule}
  Il existe un unique morphisme d'opérades $\kappa$ de $\arbu$ dans
  $\mld$ tel que
  \begin{align}
    \kappa( [2 \tri 1] ) & =  \frac{1}{u_1 (u_1+u_2)},\\
    \kappa( [1][2] ) & = \frac{1}{u_1 u_2}.
  \end{align}
\end{proposition}

\smallskip
\begin{lemma}
  \label{kappa_2}
  Soient $Q$ et $R$ deux arbustes. On a $\kappa(Q \, R)=\kappa(Q)\kappa(R)$ et 
  \begin{equation}
    \kappa(Q \tri R)=\kappa(Q)\kappa(R) \frac{\sum_{i \in Q} u_i}{\sum_{i \in Q\sqcup R} u_i}.
  \end{equation}
\end{lemma}
\begin{preuve}
  Ceci résulte du fait que $\kappa$ est un morphisme d'opérade.
\end{preuve}

On va donner maintenant une description directe du morphisme $\kappa$.
On a besoin pour cela de quelques notations.

Deux sommets ramifiés $i$ et $i'$ d'un arbuste $P$ sont dits
\textbf{équivalents} si les ensembles des buts des arêtes sortantes de
$i$ et de $i'$ sont identiques. On note $\ram(P)$ l'ensemble des
classes d'équivalence de sommets ramifiés d'un arbuste $P$. Dans la
figure \ref{arbuste5}, on a colorié les sommets ramifiés avec un
couleur commune par classe d'équivalence. Si $r\in\ram(P)$ est une
classe d'équivalence, on note $r^-$ l'ensemble commun des buts des
arêtes sortantes de $r$.

Une partie $J$ de $P$ est un \textbf{idéal supérieur} si le
complémentaire de $J$ dans $P$ est un arbuste, pour la restriction de
la fonction hauteur.

Soit $S$ une partie d'un arbuste $P$. Il existe un idéal supérieur
minimal contenant $S$, noté $\ideal{P}{S}$ et nommé l'idéal supérieur
engendré par $S$ dans $P$. C'est la partie de $I$ formée des éléments
$j$ tels que : tout chemin dans $P$ qui descend de $j$ vers un sommet
de hauteur $0$ passe par un élément de $S$. Ceci revient à dire que,
pour tout ordre total $\leq$ compatible avec $P$, il existe $i$ dans
$S$ tel que $j \geq i$.

Pour une partie $S$ de $I$, notons $\som{S}$ la somme des $u_k$ pour
$k \in S$.

Soit $P$ un arbuste. Définissons une fraction $f_P$ par la formule
\begin{equation}
  f_P = \frac{1}{\prod_{i \in I} \som{\ideal{P}{\{i\}}}}
\prod_{r \in \ram(P)} \frac{ \som{\ideal{P\setminus\ideal{P}{r}}{r^-} }}
{  \som{\ideal{P}{r^-}}}.
\end{equation}

\begin{proposition}
  Pour tout arbuste $P$, la fraction $f_P$ est réduite, sans
  multiplicité, et produit de facteurs linéaires.
\end{proposition}
\begin{preuve}
  En effet, si on compare les facteurs de $f_P$ selon l'ensemble des
  sommets de hauteur minimale qu'ils contiennent (qui est soit un
  singleton $i$ soit de la forme $r^-$ avec $r\in\ram(P)$), les seules
  coïncidences possibles sont entre les deux facteurs (un dans le
  numérateur et un dans le dénominateur) associés à un même
  $r\in\ram(P)$, et ces deux facteurs sont distincts, car chaque
  élément de $\ideal{P}{r}$ est dans l'une et pas dans l'autre.
\end{preuve}

Par exemple, en abrégeant $u_A$ en $A$, $u_B$ en $B$, etc, la fraction
associée à l'arbuste de droite de la figure \ref{composi} est la
suivante :
\begin{equation}
  \frac{(B+E+F+G)(F+G)}{A\,B\,E\,F\,G\,(A+B+C+E+F+G) (A+B+E+F+G)(E+F+G)}.
\end{equation}
Dans ce cas, on a deux classes de sommets ramifiés : $\{A\}$ et $\{E\}$.

\begin{proposition}
  \label{formule_moule}
  Le morphisme $\kappa$ de $\arbu$ dans $\mld$ est donné par $\kappa(P)=f_P$.
\end{proposition}

\begin{preuve}
  Par récurrence. L'énoncé est vrai en degré $1$. Soit donc $P$ un arbuste
  à $n$ sommets avec $n\geq 2$.

  Supposons d'abord $P$ non connexe. Alors $P$ s'écrit comme union
  disjointe d'arbres connexes. Chaque $r$ dans $\ram(P)$ est contenu
  dans un composante connexe, l'ensemble $r^-$ est dans la même
  composante, et l'idéal supérieur engendré par une partie contenue dans une
  composante connexe est contenu dans cette composante connexe. La
  formule pour $f_P$ est donc clairement le produit des formules sur
  les composantes connexes. Comme $\kappa$ est aussi multiplicatif par
  le Lemme \ref{kappa_2}, on en déduit par hypothèse de récurrence le
  résultat voulu pour $P$.

  Supposons maintenant $P$ connexe. Alors $P$ s'écrit $Q \tri R$, pour
  deux arbustes $Q$ et $R$. Si $P$ a un seul sommet de hauteur nulle,
  on peut et on va supposer que $Q$ est réduit à ce sommet.

  Considérons d'abord l'idéal supérieur $\ideal{P}{\{j\}}$ engendré
  par un singleton $\{j\}$ dans $P$. Si $j$ est dans $R$, alors
  $\ideal{P}{\{j\}}$ s'identifie à $\ideal{R}{\{j\}}$, car tous les
  sommets au dessus de $j$ sont dans $R$. Si $j$ est dans $Q$, on
  distingue deux cas. Si $j$ est le seul sommet de hauteur nulle de
  $P$, alors $Q$ est un singleton, $\ideal{P}{\{j\}}=P$, et
  $\ideal{Q}{\{j\}}=Q=\{j\}$. Sinon, il existe $i\not=j $ de hauteur
  nulle dans $P$, ce qui entraîne pour tout sommet de $R$ l'existence
  d'un chemin vers $i$ et évitant $j$. Donc aucun sommet de $R$ ne
  peut appartenir à $\ideal{P}{\{j\}}$. Par conséquent, on a
  $\ideal{P}{\{j\}}=\ideal{Q}{\{j\}}$.

  Considérons ensuite les facteurs de $f_P$ associés aux sommets
  ramifiés. Il y a une inclusion de $\ram(Q)\sqcup \ram(R) \subset
  \ram(P)$. On voit aisément que, pour chaque $r$ dans $\ram(R)\subset
  \ram(P)$, la contribution de $r$ à la fraction $f_P$ est identique à
  la contribution de $r$ à la fraction $f_R$.

  On va distinguer selon qu'il y a un ou plusieurs sommets de hauteur
  $0$ dans $P$.
 
  Supposons d'abord qu'il y a un seul sommet de hauteur $0$ dans $P$.
  Alors $Q$ est un singleton et on a $\ram(P)=\ram(R)$.

  Supposons maintenant qu'il y a plusieurs sommets de hauteur $0$ dans
  $P$. Dans ce cas, on a $\ram(P)=\ram(Q) \sqcup \ram(R) \sqcup
  \{r_0\}$, où $r_0$ est l'ensemble des sommets de hauteur $0$ dans
  $R$. On voit aisément que, pour chaque $r$ dans $\ram(Q)\subset
  \ram(P)$, la contribution de $r$ à la fraction $f_P$ est identique à
  la contribution de $r$ à la fraction $f_Q$. La partie $r_0^-$ est
  l'ensemble des sommets de hauteur nulle de $P$. L'idéal supérieur
  engendré par $r^-_0$ est donc $P$. Le facteur associé à $r_0$
  dans $f_P$ est $\frac{\som{Q}}{\som{P}}$. 

  Dans tous les cas considérés, on obtient donc que
  \begin{equation}
    f_{Q \tri R}=f_Q f_R \frac{\som{Q}}{\som{P}}.
  \end{equation}
  ce qui démontre que $f_P=\kappa(P)$, en utilisant le Lemme \ref{kappa_2}.
\end{preuve}

\subsection{Injectivité des morphismes}

\paragraph{Les composantes connexes des fractions}

On montre dans ce paragraphe qu'on peut reconstituer la partition de
$I$ en composantes connexes de $P$ à partir de la fraction $f_P$.

\begin{lemma}
  Soit $P$ un arbuste connexe. Le dénominateur de $f_P$ contient le
  facteur $\sum_{i\in I} u_i$.
\end{lemma}
\begin{preuve}
  Si $P$ a un seul sommet $i$ de hauteur nulle, l'idéal supérieur
  engendré par $i$ est $I$. Le terme correspondant au sommet $i$ dans le
  dénominateur est donc $\sum_{i\in I} u_i$.

  Si $P$ a plusieurs sommets de hauteur nulle, il existe des sommets
  de hauteur $1$ couvrant tous ces sommets par la proposition
  \ref{existe_couvrant}. Ces sommets sont ramifiés et équivalents. Le
  terme correspondant à cette classe d'équivalence de sommets ramifiés
  dans le dénominateur est $\sum_{i\in I} u_i$.
\end{preuve}

Soit $f$ une fraction en les variables $(u_i)_{i\in I}$. On définit
une partition $\pi_f$ de $I$ comme suit : $i$ et $j$ sont dans la même
part de $\pi_f$ si il existe un facteur du dénominateur de $f$ qui
fait intervenir à la fois $u_i$ et $u_j$.

\begin{lemma}
  \label{c_connexe}
  Soit $P$ un arbuste sur $I$. Les composantes connexes de $P$ sont les
  parts de $\pi_{f_P}$.
\end{lemma}

\begin{preuve}
  Soit $\pi_P$ la partition de $I$ formée par les composantes de $P$.
  Si $P$ n'est pas connexe, la fraction $f_P$ est le produit des
  fractions associées aux composantes connexes de $P$. En particulier,
  le dénominateur se factorise selon les parts de $\pi_P$. Donc
  $\pi_{f_P}$ est un raffinement de $\pi_P$. Réciproquement, par le
  lemme précédent, $\pi_P$ est un raffinement de $\pi_{f_P}$.
\end{preuve}

\paragraph{De la fraction à la fonction hauteur}

On montre dans ce paragraphe qu'on peut reconstruire la partition de
$I$ selon la hauteur à partir de la fraction $f_P$.

On remarque d'abord que, comme l'opérade $\zinb$ s'injecte dans
$\mld$, on peut reconstituer à partir de la fraction $f_P$ l'ensemble
des ordres totaux compatibles avec $P$. On va utiliser cet ensemble
pour reconstruire la partition selon la hauteur.

\begin{lemma}
  \label{niveau0}
  Les éléments de $I$ de hauteur $0$ dans $P$ sont exactement ceux qui peuvent
  apparaître comme minimum dans un ordre total compatible avec $P$.
\end{lemma}

\begin{preuve}
  Si $i$ n'est pas de hauteur nulle et $\leq$ est un ordre total
  compatible, alors il existe au moins un $j$ tel que $j < i$, donc
  $i$ n'est pas minimal. Réciproquement, si $i$ est de hauteur nulle,
  on peut toujours trouver un ordre compatible avec $P$ ayant $i$ pour
  minimum, par exemple en ordonnant les sommets par hauteur
  croissante, puis arbitrairement à hauteur fixée.
\end{preuve}

\begin{proposition}
  \label{prop_hauteur}
  On peut reconstruire à partir de $f_P$ la partition en hauteur de $I$
  associée à $P$.
\end{proposition}
\begin{preuve}
  On commence par retrouver l'ensemble des ordres totaux compatibles
  avec $P$. On en déduit les sommets de hauteur $0$ par le Lemme
  \ref{niveau0}. En supprimant ces sommets dans les ordres totaux, on
  obtient l'ensemble des ordres totaux compatibles avec l'arbuste
  obtenu par restriction de $P$ à la hauteur au moins $1$. On en
  déduit les sommets de hauteur $1$ par le Lemme \ref{niveau0}. On procède
  ainsi de suite par hauteur croissante, en éliminant les sommets dont
  la hauteur est déjà identifiée.
\end{preuve}

\paragraph{De la fraction à l'arbuste}

On montre enfin dans ce paragraphe qu'on peut reconstituer entièrement
l'arbuste $P$ à partir de la fraction $f_P$.

\begin{proposition}
  \label{injectif}
  Le morphisme $\gamma$ de $\arbu$ dans $\zinb$ est injectif, ainsi
  que le morphisme composé $\kappa$ de $\arbu$ dans $\mld$.
\end{proposition}

\begin{preuve}
  Il suffit de montrer que le morphisme $\kappa$ de $\arbu$ dans $\mld$
  est injectif. Il faut donc montrer comment reconstruire un arbuste
  $P$ à partir de sa fraction $f_P$. On donne un algorithme pour ceci,
  par récurrence.

  La fraction $1/u_1 \in \mld(\{1\})$ correspond évidemment au seul
  arbuste sur $\{1\}$.

  Étape 1 : on sait reconnaître les composantes connexes par le Lemme
  \ref{c_connexe}: les composantes connexes d'un arbuste $P$ sont
  exactement celles du dénominateur de la fraction $f_P$. On obtient
  par factorisation les fractions associées aux composantes connexes.

  Étape 2 : on peut donc supposer la fraction $f$ connexe. Elle
  contient donc la somme de toutes les variables en dénominateur. On
  sait aussi reconstruire la fonction hauteur, par la Proposition
  \ref{prop_hauteur}.

  On distingue ensuite deux cas.

  \textbf{Cas 1} : il y a un seul sommet $i$ de hauteur nulle.

  Dans ce cas, $P$ s'écrit comme $\{i\} \tri Q$ et la fraction $f_P$
  est égale à $f_Q \times \frac{1}{u[Q]+u_i}$ par le Lemme
  \ref{kappa_2}. On en déduit $f_Q$. On reconstitue alors $Q$ par
  récurrence et on obtient donc $P$.

  \textbf{Cas 2} : il y a plusieurs sommets de hauteur nulle.
  
  Alors il existe un élément $r\in\ram(P)$ correspondant aux élément
  de hauteur $1$ qui couvrent tous les élément de hauteur $0$. Dans
  $f_P$, le facteur $\alpha_r$ correspondant du numérateur est le seul
  qui contienne tous les éléments de hauteur nulle. Dans cette
  situation, $P$ s'écrit $Q \tri R$ où les sommets de $R$ sont
  exactement ceux qui n'apparaissent pas dans le facteur $\alpha_r$. En
  utilisant la formule
  \begin{equation}
    f_{Q \tri R}=f_Q f_R \frac{\som{Q}}{\som{P}},
  \end{equation}
  et en séparant les facteurs selon que leurs variables sont dans $R$
  ou non, on obtient les fractions $f_Q$ et $f_R$. On reconstitue $Q$ et $R$ par
  récurrence, donc $P$ aussi.
\end{preuve}

\section{Structure anticyclique}

L'opérade $\mld$ est munie de la structure plus riche d'opérade
anticyclique \cite{moules,moules2}. Ceci signifie qu'il existe sur
$\mld(\{1,\dots,n\}$ une action du groupe symétrique $\sym_{n+1}$ qui
étend l'action du groupe symétrique $\sym_n$ et satisfait à certaines
relations avec les compositions de l'opérade $\mld$. L'action du
groupe symétrique $\sym_{n+1}$ sur $\mld(\{1,\dots,n\})$ est
donnée comme suit. On introduit une variable supplémentaire $u_0$
vérifiant $u_0+u_1+\dots+u_n=0$. Le groupe $\sym_{n+1}$ agit par
permutations des variables $u_0,\dots,u_n$.  Toute fraction obtenue
ainsi peut s'écrire de façon unique comme une fraction en les
variables $u_1,\dots,u_n$, en éliminant $u_0$.

L 'opérade $\zinb$ est alors une sous-opérade anticyclique de $\mld$,
ce qui signifie que $\zinb(\{1,\dots,n\})$ est stable par l'action de
$\sym_{n+1}$.

Pour pouvoir parler d'opérade anticyclique ensembliste, on a besoin de
signes, car les axiomes font intervenir des signes. On définit donc
une opérade $\arbu \times \ZZ_{/2}$ en introduisant un signe $+1$ ou
$-1$ pour chaque arbuste et en prolongeant les compositions en
multipliant les signes.

On peut alors prolonger de manière unique les morphismes injectifs
d'opérades ensemblistes $\gamma$ et $\kappa$ de $\arbu$ dans $\zinb$
et $\mld$ en morphismes de $\arbu \times \ZZ_{/2}$ dans $\zinb$ et
$\mld$ qui envoient $-x$ sur $-\gamma(x)$ et $-\kappa(x)$ pour tout
$x$. On vérifie facilement que ces morphismes restent injectifs.

\begin{proposition}
  La sous opérade ensembliste $\arbu\times \ZZ_{/2}$ est une sous-opérade
  ensembliste anticyclique de $\zinb$ et $\mld$.
\end{proposition}

\begin{preuve}
  Comme $\arbu \times \ZZ_{/2}$ est engendrée en degré $2$, il suffit
  de vérifier que les images des générateurs de $\arbu\times \ZZ_{/2}$
  dans $\mld$ sont bien stables sous l'action du groupe symétrique
  $\sym_3$, ce qui est facile.
\end{preuve}

Par conséquent, il existe une action du groupe symétrique $\sym_{n+1}$
sur les arbustes signés sur l'ensemble $\{1,\dots,n\}$, pour tout
$n\geq 1$. On appellera cette action l'\textbf{action anticyclique}.

Dans la figure \ref{arbuste5}, on a groupé ensemble les classes
d'isomorphisme (\textit{i.e.} les orbites pour l'action de $\sym_n$)
d'arbustes connexes qui sont (au signe près) dans la même orbite pour
$\sym_{n+1}$. Dans chacune de ces orbites pour $\sym_{n+1}$, il peut y
avoir aussi des arbustes non connexes, qui ne sont pas représentés
dans la figure \ref{arbuste5}.

\begin{proposition}
  \label{preserve_ramif}
  L'action anticyclique sur les arbustes signés préserve le nombre de
  classes d'équivalence de sommets ramifiés.
\end{proposition}
\begin{preuve}
  Par définition de son action sur les fractions, elle préserve le
  degré du numérateur de $f_P$, ce qui entraîne le résultat.
\end{preuve}

\textbf{Remarque} : pour distinguer les orbites de cette action, on dispose
d'un invariant plus fin que le degré du numérateur de $f_P$. Cet
invariant est une paire de multi-ensembles d'entiers compris entre $1$
et $(n+1)/2$. Pour cela, on spécialise les variables $u_i$ en $1$.
Chaque facteur du numérateur et du dénominateur de $f_P$ donne un
nombre. On a donc un multi-ensemble pour le numérateur et un autre
pour le dénominateur. On remplace dans ces multi-ensembles un nombre
$k$ par $n+1-k$ si il est plus grand que $(n+1)/2$. La paire de
multi-ensembles ainsi obtenue est alors invariante pour l'action
du groupe symétrique $\sym_{n+1}$.

\subsection{Action sur les forêts d'arbres enracinés}

Par la proposition \ref{preserve_ramif}, on obtient en particulier une
action du groupe symétrique $\sym_{n+1}$ sur les forêts signées
d'arbres enracinés sur l'ensemble $\{1,\dots,n\}$, car les forêts sont
exactement les arbustes sans sommets ramifiés.

Dans cette section, on va décrire cette action.

Si $F$ est une forêt formée des arbres enracinés $T_1,\dots,T_k$, on
note $B_+(j,F)$ ou $B_+(j,T_1,\dots,T_k)$ l'arbre enraciné obtenu en
greffant les racines des arbres composant $F$ sur une nouvelle racine
commune $j$. Si $T$ est un arbre enraciné, on note $B_-(T)$ la forêt
obtenue en supprimant la racine de $T$.

\begin{lemma}
  \label{echange}
  Soit $F$ une forêt formée des arbres enracinés $T_1,\dots,T_k$. Soit
  $i$ la racine de $T_1$. Alors, pour l'action anticyclique, on a
  \begin{equation}
    \tau_{0,i} F = -B_-(T_1) \sqcup B_+(i,T_2,\dots,T_k),
  \end{equation}
  où $\tau_{0,i}$ est la transposition de $0$ et $i$.
\end{lemma}

\begin{preuve}
  Dans la fraction $f_F$, le sommet $i$ apparaît seulement dans le
  facteur correspondant à $T_1$ tout entier. En échangeant $i$ et $0$,
  et en éliminant $u_0$ par la relation $\sum_{j=0}^n u_j =0$, on obtient la
  fraction
  \begin{equation}
    -f_F \times \frac{\sum_{j\in T_1} u_j}{u_i+\sum_{j \not\in T_1} u_j}.
  \end{equation}
  Cette fraction correspond (avec un signe moins) à la forêt $B_-(T_1)
  \sqcup B(i,T_2,\dots,T_k)$.  
\end{preuve}

On note l'ensemble $C(n+1)$ des arbres enracinés signés sur l'ensemble
$\{0,1,\dots,n\}$ modulo la relation suivante : Si $(T,r)$ est un
arbre enraciné en $r$ et $r'$ un sommet adjacent à $r$ alors
\begin{equation}
  \label{negaswap}
  (T,r)=-(T,r').
\end{equation}
Le groupe symétrique $\sym_{n+1}$ agit naturellement sur $C(n+1)$. On
remarque, en utilisant la simple connexité des arbres, que chaque
élément de $C(n+1)$ admet un unique représentant enraciné en $0$. Le
cardinal de $C(n+1)$ est donc $2 (n+1)^{n-1}$.

On définit alors une application $B_0$ de l'ensemble des forêts
signées dans $C(n+1)$ comme suit : $B_0(F)=B_+(0,F)$ et
$B_0(-F)=-B_0(F)$. C'est une bijection, qui respecte les actions du
groupe $\sym_n$.

\begin{proposition}
  \label{action_forets}
  L'application $B_0$ est compatible avec les actions de $\sym_{n+1}$
  sur les forêts signées d'arbres enracinés et sur $C(n+1)$.
\end{proposition}

\begin{preuve}
  Il suffit de considérer les transpositions $\tau_{0,i}$ avec
  $i\in\{1,\dots,n\}$. On procède par récurrence sur la distance entre
  $0$ et $i$ dans l'arbre $B(0,F)$.

  Si cette distance vaut $1$, $i$ est la racine d'un des arbres
  composant $F$. Si on fait agir la transposition $\tau_{0,i}$ sur
  $B(0,F)$, on obtient un arbre enraciné en le sommet $i$, voisin du
  sommet $0$. Par la relation (\ref{negaswap}), on change sa racine
  pour la ramener en $0$ en introduisant un signe moins. Le résultat
  est
  \begin{equation}
    -B(0,B_-(T_1)
    \sqcup B(i,T_2,\dots,T_k)),
  \end{equation}
  qui est bien $B_0(\tau_{0,i}(F)$ par le lemme \ref{echange}.

  Si la distance entre $0$ et $i$ est au moins $2$, soit
  $0,i_1,i_2,\dots,i_\ell,i$ l'unique chemin entre $0$ et $i$ dans
  $B(0,F)$. Dans ce cas, on a
  \begin{equation}
    \tau_{0,i}= \tau_{i,i_1}\tau_{0,i}\tau_{0,i_1}.
  \end{equation}

  Par le cas précédent (distance $1$), on a donc
  \begin{equation}
    \tau_{0,i} B_0(F)=\tau_{i,i_1} \tau_{0,i} B_0(\tau_{0,i_1} F).
  \end{equation}
  Mais la distance dans $\tau_{0,i_1} F$ entre $0$ et $i$ est moindre
  que dans $F$, donc par hypothèse de récurrence, on a
  \begin{equation}
    \tau_{0,i} B_0(F)=\tau_{i,i_1} B_0(\tau_{0,i} \tau_{0,i_1} F).
  \end{equation}
  Comme $B_0$ est compatible avec l'action de $\sym_{n}$, on en déduit
  le résultat voulu.
\end{preuve}


\subsection{Une déformation du morphisme de $\arbu$ dans $\mld$}

A titre de remarque, on décrit ici des morphismes d'opérade
anticyclique des arbustes dans les moules, qui généralisent celui
utilisé dans le reste de l'article.

Soit $\tfun(u)$ une fonction rationnelle d'une variable $u$.
\begin{proposition}
  La fonction
  \begin{equation}
    \frac{\tfun(u_1)\tfun(u_2)}{u_1 u_2 \tfun(u_1+u_2)}
  \end{equation}
  définit un élément commutatif et associatif dans $\mld(2)$, donc
  un morphisme $\kappa_{\tfun}$ de $\comm$ dans $\mld$.
\end{proposition}
\begin{preuve}
  La preuve est une simple vérification de l'associativité.
\end{preuve}

\begin{proposition}
  Il existe une unique extension du morphisme $\kappa_{\tfun}$ en un
  morphisme d'opérade anticyclique $\kappa_{\tfun}$ de $\arbu$ dans
  $\mld$.
\end{proposition}
\begin{preuve}
  Pour montrer cela, il suffit de vérifier que les images des
  générateurs de $\arbu$, qui sont imposées par la structure
  anticyclique, vérifient la relation (\ref{relaNAP}).
\end{preuve}

Le cas $\tfun(u)=1$ est celui étudié dans le reste de l'article.

\section{L'opérade duale}

En partant de la présentation quadratique binaire de l'opérade $\arbu$
donnée par le théorème \ref{presentation}, on peut calculer (par la
procédure habituelle, voir par exemple \cite[Appendix
B]{loday_lecture_notes}) une présentation de l'opérade quadratique
binaire duale $\arbu^{!}$. On va décrire les algèbres sur cette
opérade, ce qui équivaut à décrire la présentation.

Une algèbre sur l'opérade duale $\arbu^{!}$ est donnée par un
crochet $[\,,\,]$ antisymétrique et un produit $*$
(non commutatif), vérifiant les axiomes suivants :
\begin{align}
  \label{jacobi}
  [a,[b,c]]+  [b,[c,a]]+  [c,[a,b]] & =0,\\
  [a,b]*c&=0,\\
  a*(b*c)&=0,\\
  [a*b,c]&=0,\\
  \label{loi_module}
  a*[b,c]&=(a*b)*c-(a*c)*b.
\end{align}

On remarque que le crochet est un crochet de Lie par (\ref{jacobi}),
et que le produit $*$ n'est pas associatif, mais définit un module à
droite pour le crochet de Lie par (\ref{loi_module}).

Les $\arbu^{!}$-algèbres libres ont une description simple. En termes
d'opérades, ceci s'exprime comme suit.

\begin{proposition}
  On a
  \begin{equation}
    \arbu^{!}(n)=
    \begin{cases}
      \ass(1)&\text{ si }n=1,\\
      \lie(n)\oplus\ass(n)&\text{ si }n\geq 2,
    \end{cases}
  \end{equation}
  où $\lie(n)=\lie(\{1,\dots,n\})$ et $\ass(n)=\ass(\{1,\dots,n\})$
  sont les composantes des opérades classiques $\lie$ et $\ass$.
\end{proposition}

\begin{preuve}
  C'est une conséquence facile de la présentation.
\end{preuve}

\medskip

L'opérade $\arbu^{!}$ hérite par dualité d'une structure anticyclique
provenant de celle de $\arbu$. Ceci donne une notion de forme
antisymétrique invariante $ \langle \,,\, \rangle$ pour les algèbres sur
$\arbu^{!}$. La condition d'invariance est la suivante :
\begin{equation}
  \langle a*b,c \rangle = - \langle [b,c] , a \rangle.
\end{equation}

\bibliographystyle{alpha}
\bibliography{cloture}

\end{document}